\documentclass[12pt,fleqn,letterpaper]{article}
\usepackage{a4wide,amsfonts,amsmath,latexsym,amssymb,euscript,graphicx,units,mathrsfs}

\usepackage{pstricks, pst-plot, pst-node, pst-tree}
\usepackage{graphicx}
\usepackage{color}
 \usepackage{dsfont}
\usepackage{geometry}
\usepackage{natbib}
\usepackage{amssymb}
\usepackage{amssymb}
\usepackage[T1]{fontenc}
\usepackage{latexsym}
\usepackage{xypic}
  \usepackage{placeins}
\usepackage{eufrak}
\usepackage{euscript}
\usepackage{amsfonts,amsmath}
\usepackage{verbatim}
\usepackage{fancyhdr}
\usepackage[english]{babel}
\usepackage{units}
\usepackage[colorlinks=true,linkcolor=red,citecolor=blue]{hyperref}

\newtheorem{theorem}{Theorem}[section]

\newtheorem{remark}[theorem]{Remark}

\numberwithin{equation}{section} 
\begin{document}

\title{New Entropy Estimator with an Application to Test of Normality}

\author{Salim BOUZEBDA$^{1}$,\footnote{e-mail: salim.bouzebda@upmc.fr~~(corresponding author)} \hbox{ } Issam
ELHATTAB$^{1}$,\footnote{e-mail:
issam.elhattab@upmc.fr}\\ Amor KEZIOU$^{1}$ \footnote{e-mail:
amor.keziou@upmc.fr} \hbox{ } and
Tewfik LOUNIS $^2$\footnote{e-mail: tewfik.lounis@gmail.com}
\\
L.S.T.A.,  Universit\'e Pierre et Marie Curie\\
    4 place Jussieu
    75252 Paris Cedex 05\\Laboratoire de Math\'ematiques Nicolas Oresme\\ CNRS UMR 6139
Université de Caen$^2$ }
\date{}
\maketitle

\begin{abstract}
\noindent In the present paper we propose a new estimator of entropy
based on smooth estimators of quantile density. The consistency and
asymptotic distribution of the proposed estimates are obtained. As a consequence,
a new test of normality is proposed. A small power comparison is provided. 
A simulation study for the comparison, in terms of mean squared error,
of all estimators under study is
performed.
\end{abstract}

\section{Introduction and estimation}\label{6-1mai}
\noindent Let $X$ be a random variable [r.v.] with cumulative
distribution function [cdf] $F( x):=\mathbb{P}(X\leq  x)~~\mbox{ for } ~~
x\in\mathbb{R}$ and a density function $f(\cdot)$ with respect to
Lebesgue measure on $\mathbb{R}$. Then its differential  (or
Shannon) entropy is defined by
\begin{equation}\label{6-entropy}
H(X) := - \int_{\mathbb{R}} f (x) \log f (x)dx,
\end{equation}
where $dx$ denotes Lebesgue measure on $\mathbb{R}$.
We assume that $H(X)$ is
properly defined by the integral (\ref{6-entropy}),
 in the sense that
\begin{equation}\label{6-entropyfinite}
\vert H(X)\vert<\infty.
\end{equation}
The concept of differential entropy was introduced in Shannon's
original paper [\cite{shannon1948}]. Since this early epoch, the
notion of entropy has been the subject of great theoretical and
applied interest. Entropy concepts and principles play a
fundamental  role in many applications, such as statistical
communication theory [\cite{Gallager1968}], quantization theory
[\cite{Reny1959}], statistical decision theory
[\cite{kullback1959}], and  contingency table analysis
[\cite{Gokhale1978}]. \cite{Csizar1962} introduced the concept of
convergence in entropy and showed that the latter convergence
concept implies convergence in $\mathcal{L}_1$. This property
indicates that entropy is a useful concept to measure
``\emph{closeness in distribution}'', and also justifies
heuristically the usage of sample entropy as test statistics when
designing entropy-based tests of goodness-of-fit. This line of
research has been pursued by \cite{vasicek1976},
\cite{Prescott1976}, \cite{DudewiczMeulen1981},
\cite{Gokhale1983},
 \cite{Ebrahimi1992} and \cite{Esteban2001} [including the
references therein].  The idea here is that many families of
distributions are characterized by maximization of entropy subject
to constraints (see, e.g., \cite{Jaynes1957} and
\cite{Verdugo1978}). There is a huge literature on the Shannon's entropy and its applications.
It is not the purpose of this paper to survey this extensive literature. We can refer
 to \cite{coverthomas2006} (see their Chapter
8), for a comprehensive overview of differential entropy and their
mathematical properties.\vskip5pt\noindent In the literature, Several proposals have
been made to estimate entropy. \cite{dmitrievtarasenko1973} and
\cite{ahmadlin1976} proposed estimators of the entropy using
kernel-type estimators of the density $f(\cdot)$.
\cite{vasicek1976} proposed an entropy estimator based on
spacings. Inspired by the work of \cite{vasicek1976}, some authors
[\cite{vanes1992}, \cite{Correa1995} and \cite{Wieczorkowski1999}]
proposed modified entropy estimators, improving in some respects
the properties of Vasicek's estimator (see Section
\ref{comparison} below). The reader finds in
\cite{beirlantdudewiczgyorfivandermeulen1997} detailed accounts of
the theory as well as surveys for entropy
estimators.\vskip5pt\noindent This paper aims to introduce a new
entropy estimator and obtains its asymptotic properties.
Comparison Simulations indicate that our estimator produces
smaller mean squared error than the other well-known competitors
considered in this work.  As a consequence, we propose a new test
of normality.
\noindent First,
we introduce some definitions and notations.
For each distribution function $F(\cdot)$, we define the
quantile function by
\begin{equation*}
Q(t):=\inf \{x: F(x)\geq t\},~~ 0<t<1.
\end{equation*}
Let $$x_F:=\sup\{x: F(x)=0\}~~\mbox{and}~~x^F:=\inf\{x: F(x)=1\},~~
-\infty\leq x_F<x^F\leq \infty.$$  We assume that the distribution
function $F(\cdot)$ has a density $f(\cdot)$ (with respect to
Lebesgue measure on $\mathbb{R}$), and that $f(x)>0$ for all
$x\in(x_F,x^F)$. Let $$q(x):=dQ(x)/dx=1/f(Q(x)),~~ 0<x<1,$$ be the
quantile density function [qdf]. The entropy $H(X)$, defined by
(\ref{6-entropy}), can be expressed in the form of
quantile-density functional as
\begin{equation}\label{6-entropyq}
H(X)=\int_{[0,1]}\log\left(\frac{d}{dx}Q(x)\right) dx=\int_{[0,1]}\log\big(q(x)\big) dx.
\end{equation}
The Vasicek's estimator was constructed by replacing the quantile
function $Q(\cdot)$ in (\ref{6-entropyq}) by the empirical quantile function
and using a difference operator instead of the differential one. The derivative
$(d/dx)Q(x)$ is then estimated
by a function of spacings.\vskip5pt\noindent In this work, we construct our estimator of entropy by replacing $q(\cdot)$,
in (\ref{6-entropyq}), by an appropriate estimator $\widehat{q}_{n}(\cdot)$ of $q(\cdot)$.
We shall consider the kernel-type estimator of $q(\cdot)$ introduced by \cite{FalkM1986}
and studied by \cite{chengParzen1997}.
Our choice is motivated by the well asymptotic behavior properties of this estimator.
\cite{chengParzen1997} were established the asymptotic properties of $\widehat{q}_{n}(\cdot)$
on all compact $U\subset ]0,1[$, which avoids the boundary problems. Since the entropy is definite as an integral on $]0,1[$ of a
 functional of $q(\cdot)$, it is not suitable to substitute directly  $\widehat{q}_{n}(\cdot)$ in
 (\ref{6-entropyq}) to estimate $H(X)$. To circumvent the boundary effects, we will proceed as follows.
\noindent We set for small $\varepsilon\in]0,1/2[$, $$H_\varepsilon(X):=\varepsilon\log\big
(q(\varepsilon)\big)+\varepsilon\log\big(q(1-\varepsilon)\big)+
    \int_{\varepsilon}^{1-\varepsilon}\log\big(q(x)\big) dx.$$
In view (\ref{6-entropyq}), we can see that
\begin{equation}\label{Keyidea}
    | H(X)-H_\varepsilon(X)|=o\big(\eta(\varepsilon)\big),
\end{equation}
where $\eta(\varepsilon)\rightarrow0,$ as $\varepsilon\rightarrow0$. The choice of $\varepsilon$
close to zero guaranteed the closeness of $H_\varepsilon(X)$ to $H(X)$, then the problem of the estimation of $H(X)$
is reduced to estimate $H_\varepsilon(X)$.

\noindent Given an independent and
identically distributed random [i.i.d.] sample $X_1,\dots,X_n$, and let $\varepsilon\in]0,1/2[$,
an estimator of $H_\varepsilon(X)$ can be defined as
\begin{equation}\label{6-entropy-estimator}\widehat{H}_{\varepsilon;n}(X)=
\varepsilon\log\left(\widehat{q}_n\left(\varepsilon\right)\right)
+\varepsilon\log\left(\widehat{q}_n\left(1-\varepsilon\right)\right)
+\int_{\varepsilon}^{1-\varepsilon}\log\big(\widehat{q}_n(x)\big) dx,
\end{equation}
\noindent  where the
estimator $\widehat{q}_{n}(\cdot)$ of the qdf $q(\cdot)$ is defined as follows. Let  $X_{1;n}\leq \cdots\leq X_{n;n} $ denote
the order statistics of $X_1,\ldots,X_n$. The empirical quantile function
$Q_n(\cdot)$ based upon these random variables is given by
\begin{equation}
Q_n(t):=X_{k;n}, ~~(k-1)/n<t\leq k/n,~~k=1,\ldots,n.
\end{equation}
Let $\{K_n(t,x), (t,x)\in[0,1]\times]0,1[\}_{n\geq 1}$ be a sequence of kernels
and $\{\mu_n(\cdot)\}_{n\geq1}$ a sequence of $\sigma$-finite measures on $[0,1]$.
A smoothed version of $Q_n(\cdot)$ (see, e.g., \cite{chengParzen1997}) can be defined as
\begin{equation*}
    \widehat{Q}_n(t):=\int_0^1Q_n(x)K_n\left(t,x\right)d\mu_n(x),
      \quad t\in]0,1[.
\end{equation*}
Finally, we estimate $q(\cdot)$ by
\begin{eqnarray}\label{6-generalformqdf}
    \widehat{q}_n(t):=\frac{d}{dt}\widehat{Q}_n(t)
    =\frac{d}{dt}\int_0^1Q_n(x)K_n\left(t,x\right)d\mu_n(x),
    \quad t\in]0,1[.
\end{eqnarray}
Clearly, in order to obtain a meaningful qdf estimator in this way, the sequence
of kernels $K_n\left(\cdot,\cdot\right)$ must satisfy certain differentiability
conditions, and together with the sequence of  measures $\mu_n(\cdot)$, must satisfy
certain variational conditions. These conditions will be detailed in the next section.
A familiar example is the convolution-kernel estimator
\begin{eqnarray*}
    \widehat{q}_n(t):=\frac{d}{dt}\int_0^1h_n^{-1}Q_n(x)K\left(\frac{t-x}{h_n}\right)dx,
    \quad t\in]0,1[,
\end{eqnarray*}
where $K(\cdot)$ denotes a kernel function,
namely a measurable function integrating to 1 on $\mathbb{R}$, and has bounded derivative,
and $\{h_n\}_{n\geq1}$ is a sequence of positive
reals fulfills $h_n\rightarrow 0 $ and $nh_n\rightarrow \infty$ as $n\rightarrow\infty$.
\noindent In this case, \cite{csorgorevesz1984} and
\cite{csorgodeuvelshorvath1991} define
\begin{eqnarray*}
\bar{q}_n(t)&:=&h_n^{-1}\int_{1/(n+1)}^{n/(n+1)}K\left(\frac{t-x}{h_n}\right)~dQ_n(x)\nonumber\\
&=&h_n^{-1}\sum_{i=1}^{n-1}K\left(\frac{t-i/n}{h_n}\right)(X_{i+1;n}-X_{i;n}),~~t\in[0,1].
\end{eqnarray*}
Calculations using summation by parts show that, for all
$t\in]0,1[$,
\begin{equation*}
    \widehat{q}_n(t)=\bar{q}_n(t)+h_n^{-1}\left[K\left(\frac{t-1}{h_n}\right)X_{n;n}
            -K\left(\frac{t}{h_n}\right)X_{1;n}\right].
\end{equation*}
In the sequel, we shall consider the general family of qdf
estimators defined in (\ref{6-generalformqdf}).\\ \vskip5pt

\noindent The remainder of the present article is organized as follows.
The consistency and normality of our estimator
are discussed in the next section. Our arguments, used to establish the asymptotic normality of our estimator,
make use of an
original application of the invariance principle for the empirical quantile process.
In Section \ref{parzen}, we discuss briefly the smoothed version of Parzen's entropy estimator.
In Section \ref{comparison}, we investigate the
finite-sample performance of the newly proposed estimator and
compare the latter with the performances of existing estimators. In section \ref{simulation},
a new test of normality is proposed and compared with other competitor.
Some
concluding remarks and possible future developments are mentioned
in Section \ref{conclusion}.
To avoid interrupting the flow of the
presentation, all mathematical developments are relegated to Section
\ref{6-demo}.

\section{Main results}\label{6-mai}
\noindent Throughout $U(\varepsilon):=[\varepsilon,1-\varepsilon]$,
for an $\varepsilon\in]0,1/2[$ is arbitrarily fixed (free of the sample size $n$).
The following conditions are used to establish
the main results. We recall the notations of Section \ref{6-1mai}.
\begin{enumerate}
\item[(Q.1)] The quantile density function $q(\cdot)$ is twice
    differentiable in $]0,1[$. Moreover,
    $$0<\min\{q(0),q(1)\}\leq\infty;$$
\item[(Q.2)] There exists  a constant $\varsigma>0$ such that
    $$\sup_{t\in]0,1[} \big\{t(1-t)~\vert
    J(t)\vert\big\}\leq\varsigma,$$ where
    $J(t):=d\log\big\{q(t)\big\}/dt$
    is the score function;
\item[(Q.3)] Either $q(0)<\infty$ or $q(\cdot)$ is
    nonincreasing in some interval $(0,t_*)$, and either
    $q(1)<\infty$ or $q(\cdot)$ is nondecreasing in some
    interval $(t^*,1)$, where $0<t_*<t^*<1$.
\end{enumerate}
\noindent We will make the following assumptions on the sequence of kernels
$K_{n}(\cdot,\cdot)$.
\begin{enumerate}
\item[(K.1)] For each $n\geq1$ and each $(t,x)\in U(\varepsilon)\times]0,1[$, $K_{n}(t,x)\geq0$,
and for each $t\in U(\varepsilon)$, $\int_0^1K_{n}(t,x)d\mu_n(x)=1$;
\item[(K.2)] There is a sequence $\delta_n\downarrow0$ such
    that
$$R_n:=\sup_{t\in U(\varepsilon)}\left[1-\int_{t-\delta_n}^{t+\delta_n}
K_n\left(t,x\right)d\mu_n(x)\right]\rightarrow0, \quad\mbox{as~~} n\rightarrow\infty;$$
\item[(K.3)] For any function $g(\cdot)$, that is at least three times
    differentiable in $]0,1[$, $$\widehat{g}_n(t):=\int_0^1g(x)K_n\left(t,x\right)d\mu_n(x),$$
    is differentiable in $t$ on $U(\varepsilon)$, and
        $$\sup_{t\in U(\varepsilon)}\left\vert g(t)-\int_0^1g(x)
        K_n\left(t,x\right)d\mu_n(x)\right\vert=O\left(n^{-\alpha}\right), \qquad \alpha>0;$$
        $$\sup_{t\in U(\varepsilon)}\left\vert g'(t)-\frac{d}{dt}\int_0^1g(x)
        K_n\left(t,x\right)d\mu_n(x) \right\vert=O\left(n^{-\beta}\right), \qquad \beta>0.$$
\end{enumerate}
Note that the conditions (Q.1-2-3) and (K.1-2-3), which we will
use to establish the convergence in probability of $\widehat{H}_{\varepsilon;n}(X)$, match
those used by \cite{chengParzen1997} to establish the asymptotic properties
of $\widehat{q}_n(\cdot)$.\vskip5pt\noindent
Our result concerning  consistency of $\widehat{H}_{\varepsilon;n}(X)$ defined in
(\ref{6-entropy-estimator}), where the qdf estimator is given in (\ref{6-generalformqdf}), is as follows.
\begin{theorem}\label{6-theorem1}
Let $X_1,\ldots,X_n$ be i.i.d. r.v.'s with a quantile density
function $q(\cdot)$ fulfilling {\rm (Q.1-2-3)}. Let $K_n(\cdot,\cdot)$
fulfill {\rm (K.1-2-3)}. Then, we have,
\begin{equation}\label{6-as}
|\widehat{H}_{\varepsilon;n}(X)-H(X)|=O_{\mathbb{P}}\Big(n^{-1/2}M(q;K_n)+n^{-\beta}+\eta(\varepsilon)\Big),
\end{equation}
where
\begin{eqnarray*}
  M(q;K_n) &:=& M_q\widehat{M}_n(1)\sqrt{\delta_n\log\delta_n^{-1}}+M_{q'}+\sqrt{\widehat{M}_n(q^2)R_n'(1)}+n^{-1/2}A_\gamma(n)M_q\widehat{M}_n(1),\\
  \widehat{M}_n(g) &:=& \sup_{u\in U(\varepsilon)}\int_0^1\vert g(x)K_n(u,x)\vert d\mu_n(x), \\
  R_n'(g) &:=& \sup_{u\in U(\varepsilon)}\int_{[0,1]\backslash U(\varepsilon+\delta_n)}\vert g(x)K_n(u,x)\vert d\mu_n(x), \\
  M_g &:=& \sup_{u\in U(\varepsilon)}\vert g(u)\vert,
\end{eqnarray*}
and $\delta_{n}$ is given in condition {\rm (K.2)}.
\end{theorem}%
The proof of Theorem \ref{6-theorem1} is postponed to the Section \ref{6-demo}.\\
\noindent
To establish the asymptotic normality of $\widehat{H}_{\varepsilon;n}(X)$, we will use the
following additional condition.
\begin{enumerate}
\item[(Q.4)] $\mathbb{E}\left\{\log^2\Big(q\big(F(X)\big)\Big)\right\}<\infty.$
\end{enumerate}
Let, for all $\varepsilon\in]0,1/2[$,
\begin{eqnarray*}
    \mathbb{E}\left\{\log^2\Big(q\big(F(X)\big)\Big)\right\}-
    \left\{\varepsilon\log^2\big(q(\varepsilon)\big)+\varepsilon\log^2\big(q(1-\varepsilon)\big)+
    \int_{\varepsilon}^{1-\varepsilon}\log^2\big(q(x)\big) dx\right\}=o\Big(\vartheta(\varepsilon)\Big),
\end{eqnarray*}
where $\vartheta(\varepsilon)\rightarrow0$, as $\varepsilon\rightarrow0$.
\\
The main result, concerning the normality of $\widehat{H}_{\varepsilon;n}(X)$, to be proved here may now be stated precisely as
follows.
\begin{theorem}\label{6-norm1}
Assume that the conditions  {\rm (Q.1-2-3-4)} and {\rm (K.1-2-3)} hold with $\alpha>1/2$ and $\beta>1/2$ in {\rm (K.3)}.
Then, we have,
\begin{equation}\label{6-norm}
\sqrt{n}(\widehat{H}_{\varepsilon;n}(X) -H_\varepsilon(X))-\psi_n(\varepsilon)= O_{\mathbb{P}}\left(\Big\{2\varepsilon\log\varepsilon^{-1}\Big\}^{1/2}\right)+o_{\mathbb{P}}(1),
\end{equation}
where $$\psi_n(\varepsilon):=
\int_{U(\varepsilon)}(q^{\prime}(x)/q(x))B_{n}(x)dx,$$ is a centered
Gaussian random variable with variance equal to
\begin{equation*}
    \mbox{V}ar\big(\psi_n(\varepsilon)\big)=
    \mbox{V}ar\left\{\log\Big(q\big(F(X)\big)\Big)\right\}+o\Big(\vartheta(\varepsilon)+\eta(\varepsilon)\Big).
\end{equation*}

\end{theorem}
The proof of Theorem \ref{6-norm1} is postponed to the Section \ref{6-demo}.
\begin{remark}
Condition {\rm (Q.4)} is extremely weak and is satisfied by all commonly encountered
distributions including many important heavy tailed distributions for which the moments
do not exists, the interested reader may refer to \cite{Song2000} and the references therein.
\end{remark}

\begin{remark}
We mention that the condition {\rm (K.3)} is satisfied,  when $\alpha,\beta>1/2$, for the
 difference kernels $K_n\left(t,x\right)d\mu_n(x) = h_{n}^{-1}k((t - x)/h_{n})dx$ with
 $h_{n}=O(n^{-\nu})$ where $1/4<\nu<1$. A typical example for differences kernels satisfying theses conditions is the gaussian kernel.
\end{remark}
\section{The smoothed Parzen estimator of entropy}\label{parzen}
We mention that the notations of this section are similar to that used  in \cite{Parzen1979}
 and changes have been made in order to adopt it to our setting.
Given a random sample $X_{1},\ldots,X_{n}$ of a continuous random variable $X$ with
distribution function $F(\cdot)$.
In this section we will work under the following hypothesis, for all $x\in \mathbb{R}$,
\begin{equation*}
\mathscr{H}^{*}: F(x)=F_{0}\left(\frac{x-\mu}{\sigma}\right),
\end{equation*}
where $\mu$ and $\sigma$ are parameters to be estimated, it is convenient to transform $X_{i}$
to $$\widehat{U}_{i}=F_{0}\left(\frac{X_{i}-\widehat{\mu}_{n}}{\widehat{\sigma}_{n}}\right),$$
where $\widehat{\mu}_{n}$ and $\widehat{\sigma}_{n}$ are efficient estimators of $\mu$ and $\sigma$.
 It is more convenient to base tests of the hypothesis on the deviations from identity function $t$
  of the empirical distribution function of $\widehat{U}_{1},\ldots,\widehat{U}_{n}$.
  Let denote by $\widetilde{D}_{n,0}(t)$, $0\leq t\leq 1$ the empirical quantile
  function of $\widehat{U}_{1},\ldots,\widehat{U}_{n}$; it can be expressed in terms of the
  sample quantile function $\widetilde{Q}_{n}(\cdot)$ of the original data $X_{1},\ldots,X_{n}$ by
\begin{equation*}
\widetilde{D}_{n,0}(t):=F_{0}\left(\frac{\widetilde{Q}_{n}(t)-\widehat{\mu}_{n}}{\widehat{\sigma}_{n}}\right),
\end{equation*}
where
\begin{equation*}
\widetilde{Q}_{n}(t)=n\left(\frac{i}{n}-t\right)X_{i-1;n}+n\left(t-\frac{i-1}{n}\right)X_{i;n},
 ~~\mbox{for}~~\frac{i-1}{n}\leq t\leq \frac{i}{n}, i=1,\ldots,n.
\end{equation*}
In our case, we may consider the smoothed version of the quantile density given in
\cite{Parzen1979} and defined by
$$
\widetilde{d}_{n,0}(t):=\frac{d}{d t}\widetilde{D}_{n,0}(t)=f_{0}\left(\frac{\widetilde{Q}_{n}(t)-
\widehat{\mu}_{n}}{\widehat{\sigma}_{n}}\right)\frac{d}{dt}\widetilde{Q}_{n}(t)\frac{1}{\widehat{\sigma}_{n}}.
$$
Consequently, define
$$
\widetilde{d}_{n}(t):=f_{0}(Q_{0}(t))\frac{d}{dt}\widetilde{Q}_{n}(t)\frac{1}{\widehat{\sigma}_{0,n}},
$$
where
$$
\widehat{\sigma}_{0,n}:=\int_{0}^{1}f_{0}(Q_{0}(t))\frac{d}{dt}\widetilde{Q}_{n}(t)dt.
$$
This is an alternative approach to estimating $q(\cdot)$ when one desires a goodness of fit
 test of a location scale parametric model, we refer to \cite{Parzen1990} for more details.
  We mention that in \cite[Section 7.]{Parzen1979}, the smoothed version of $\widetilde{d}_{n}(t)$ is defined, for the convolution-kernel, by
 $$
 \widetilde{d}_{n}(t):=\int_{0}^{1}\widetilde{d}_{n}(u)\frac{1}{h_{n}}K\left(\frac{t-u}{h_{n}}\right)du.
 $$
This suggests  the following estimator of entropy
\begin{equation}\label{parzen}
H_{n}(X):=\int_{0}^{1}\log(\widetilde{d}_{n}(t))dt.
\end{equation}
Under similar  conditions to that of Theorem \ref{6-theorem1}, we could derive, under $\mathscr{H}^{*}$, that, as $n\rightarrow  \infty $
\begin{equation}
|H_{n}(X)-H(X)|=o_{\mathbb{P}}(1).
\end{equation}
For more details on goodness of fit test in connection with
Shannon's entropy, the interested reader may refer to
\cite{Parzen1990}. Note that the normality  is an example of
location-scale parametric model.
\begin{remark}
Recall that there exists several ways to obtain smoothed versions of $\widetilde{d}_{n}(\cdot)$.
Indeed, we can choose the following smoothing method. Keep in mind the following definition
\begin{eqnarray*}
    \widehat{q}_n(t)=\frac{d}{dt}\widehat{Q}_n(t)
    =\frac{d}{dt}\int_0^1Q_n(x)K_n\left(t,x\right)d\mu_n(x),
    \quad t\in]0,1[.
\end{eqnarray*}
Then, we can use the following estimator
$$
\widetilde{d}_{n}^{*}(t):=f_{0}(Q_{0}(t))\widehat{q}_n(t)\frac{1}{\widehat{\sigma}_{0,n}} .
$$
Finally, we estimate the entropy under $\mathscr{H}^{*}$
\begin{equation}
H_{n}^{*}(X):=\int_{0}^{1}\log(\widetilde{d}_{n}^{*}(t))dt.
\end{equation}
Under conditions of Theorem \ref{6-theorem1}, we could derive, under $\mathscr{H}^{*}$, that, as $n\rightarrow  \infty $
\begin{equation}
|H_{n}^{*}(X)-H(X)|=o_{\mathbb{P}}(1).
\end{equation}

\end{remark}
It will be interesting to provide a complete investigation of $H_{n}(X)$ and $H_{n}^{*}(X)$
and study the test of normality of statistics based on them.  This
would go well beyond the scope of the present paper.

 \begin{remark}
 The nonparametric approach, that we use to estimate the Shannon's entropy, treats the density parameter-free and thus offers the greatest generality.
  This is in contrast with the parametric approach where one assumes that
  the underlying distribution follows a certain parametric model, and the
  problem reduces to estimating a finite number of parameters describing the model.
 \end{remark}

\section{A comparison study by simulations}\label{comparison}
\noindent Assuming that $X_1,\ldots, X_n$, is the sample, the
estimator proposed by \cite{vasicek1976} is given by
\begin{equation*}
H_{m,n}^{(V)}:=\frac{1}{n}\sum_{i=1}^n\log\left(\frac{n}{2m}\left\{X_{i+m;n}-X_{i-m;n}\right\}\right),
\end{equation*}
where $m$ is a positive integer fulfilling $m\leq n/2$. Vasicek proved
that his estimator is consistent, i.e.,
$$H_{m,n}^{(V)}\stackrel{\mathbb{P}}{\rightarrow}H(X),~~\mbox{ as }~~n\rightarrow \infty,
~~m\rightarrow\infty~~\mbox{and}~~m/n\rightarrow 0.$$ Vasicek's
estimator is also asymptotically normal under appropriate
conditions on $m$ and $n$.

\noindent \cite{vanes1992} proposed a new estimator of entropy
given by
\begin{eqnarray*}
H_{m,n}^{(Van)}&:=&-\frac{1}{n-m}\sum_{i=1}^{n-m}\log\left(\frac{n+1}{m}
\left\{X_{i+m;n}-X_{i;n}\right\}\right)\\&&+\sum_{k=m}^n\frac{1}{k}+\log(m)-\log(n+1).
\end{eqnarray*}
\cite{vanes1992} established, under suitable conditions, the consistency and
asymptotic normality of this estimator.

\noindent \cite{Correa1995} suggested a modification of Vasicek's
estimator. In order to estimate the density $f(\cdot)$ of $F(\cdot)$ in the
interval $(X_{i-m;n}, X_{i+m;n})$ he used a local linear model
based on $2m + 1$ points: $$F(X_{j;n}) = \alpha + \beta X_{j;n},
j = m-i,\dots,m+i.$$ This  produces a smaller mean squared
error estimator. The Correa's estimator is defined by
\begin{equation*}
H_{m,n}^{(Cor)}:=-\frac{1}{n}\sum_{i=1}^{n-m}\log\left(\frac{\sum_{j=i-m}^{i+m}\left\{X_{j;n}-\overline{X}_{(i)}
\right\}(j-i)}{n\sum_{j=i-m}^{i+m}\left\{X_{j;n}-\overline{X}_{(i)}\right\}}\right),
\end{equation*}
where
$$
\overline{X}_{(i)}:=\frac{1}{2m+1}\sum_{j=i-m}^{i+m}X_{j;n}.
$$
Here, $m\leq n/2$ is a positive integer, $X_{i;n} = X_{1;n}$ for
$i < 1$ and $X_{i;n} = X_{n;n}$ for $i > n$. By simulations, Correa
showed that, for some selected distributions, his estimator
produces smaller mean squared error than Vasicek's estimator and van Es' estimator. \\
\noindent \cite{Wieczorkowski1999} modified the Vasicek's estimator
by adding a bias correction. Their estimator is given by
\begin{eqnarray*}
H_{m,n}^{(WG)}&:=&H_{m,n}^{(V)}-\log n+\log(2m)-\left(1-\frac{2m}{n}\right)\Psi(2m)\\&&+\Psi(n+1)-\frac{2}{n}\sum_{i=1}^m\Psi(i+m-1),
\end{eqnarray*}
where $\Psi(x)$ is the digamma function defined by (see e.g.,
\cite{GradshteynRyzhik})
\begin{equation*}
\Psi(x):=\frac{d \log \Gamma(x)}{dx}=\frac{\Gamma^{\prime}(x)}{\Gamma^(x)}.
\end{equation*}
For integer arguments, we have 
\begin{equation*}
\Psi(k)=\sum_{i=1}^{k-1}\frac{1}{i}-\gamma,
\end{equation*}
where $\gamma=0.57721566490\ldots$ is the Euler's constant.\\
\noindent A series of experiments were conducted in order to
compare the performance of our estimator, in terms of efficiency
and robustness,  with the following entropy estimators : Vasicek's
estimator, van Es' estimator, Correa's estimator
 and Wieczorkowski-Grzegorewski's estimator.
 We provide numerical illustrations regarding the mean
squared error of each of the above estimators of the entropy
$H(X)$. The computing program codes were implemented in
\texttt{R}. We have considered the following distributions. For each distribution, we give the value of $H(X)$.
\begin{enumerate}
\item[$\bullet$] Standard normal distribution $N(0,1)$:
 $$H(X) =\log(\sigma\sqrt{2\pi  e})=\log(\sqrt{2\pi
    e}).$$
\item[$\bullet$] Uniform distribution:
     $$H(X) =0.$$
\item[$\bullet$] Weibull distribution  with the shape parameter  equal to $2$ and
the scale parameter equal to $0.5$: $$H(X) =-0.09768653.$$
    Recall that the entropy for Weibull distribution, with the shape parameter $k > 0$
    and  the scale parameter $\lambda>0$, is given by
    $$
    H(X)=\gamma\left(1-\frac{1}{k}\right)+\log\left(\frac{\lambda}{k}\right)+1.
    $$
\item[$\bullet$] Exponential distribution with parameter $\lambda=1$: $$H(X) =1-\log(\lambda)=1.$$
  \item[$\bullet$] Student's $t$-distribution  with the number of degrees of freedom equal to $1$, $3$ and $5$: 
    \begin{eqnarray*}
    H(X) &=&2.53102425, ~~~~ (\mbox{degrees of freedom equal to}~~ 1 ~~(\mbox{Cauchy distribution})),\\
    H(X) &=&1.77347757,~~~ (\mbox{degrees of freedom equal to}~~ 3), \\
    H(X) &=&1.62750267, ~~~~ (\mbox{degrees of freedom equal to}~~ 5). 
     \end{eqnarray*}
    Keeping in mind that the entropy for Student's $t$-distribution, with the number of degrees of freedom $\nu$, is given by
    $$
    H(X)=\frac{\nu+1}{2}\left(\Psi\left(\frac{1+\nu}{2}\right)-\Psi\left(\frac{\nu}{2}\right)\right)+\log
    \left(\sqrt{\nu}{\rm Beta}\left(\frac{\nu}{2},\frac{1}{2}\right)\right),
    $$
where $\Psi(\cdot)$ is digamma function.
    \end{enumerate}
For  sample sizes  $n=10$, $n=20$ and $n=50$; $5000$ samples were
generated. All the considered spacing-based estimators depend on
the parameter  $m\leq n/2$; the optimal choice of $m$ given to the
sample size $n$ is still an open problem. Here, we have chosen the
 parameter $m$ according to \cite{Correa1995}, where $m=3$ for
$n=10$ and $n=20$, and $m=4$ for $n=50$.  For our estimator
$\widehat{H}_{\epsilon,n}(X)$, we have used the standard gaussian
kernel, and the choice of the bandwidth is done by numerical
optimization of the MSE of $\widehat{H}_{\varepsilon;n}(X)$ with
respect to $h_n$. In Tables \ref{tab1}, \ref{tab2} and \ref{tab3} we have considered normal
distribution $N(0,1)$ with sample sizes $n=10$, $n=20$ and $n=50$.
In Tables \ref{tab3U}-\ref{tab3C5}, we have considered the uniform distribution,
Weibull distribution, exponential distribution with parameter $1$,
Student $t$-distribution with parameter $3$ and $5$ and Cauchy
distribution, respectively, and sample size $n=50$. In all cases,
and for each considered estimator, we compute the bias, variance
and MSE by Monte-Carlo through the 5000 replications.\vskip5pt\noindent From Tables
\ref{tab1}-\ref{tab3C5}, we can see that our estimator works better than all the
others, in the sense that the MSE is smaller.   We think that the simulation results may be substantially  ameliorated
if we choose for example the Beta or Gamma kernel,
which have the advantage to take into account the boundary
effects. It appears that
our estimator, based on the smoothed quantile density estimate,
behaves better than the traditional ones in term of efficiency.\\

\begin{table}[h!]
  \centering
  \begin{tabular}{|c|c|c|c|c|}
    \hline
     & Estimate & bias & Variance & MSE \\
    \hline
 $H_{m,n}^{(V)}$    &  0.85912656 &-0.55981198&  0.07087153&  0.38419010 \\
    \hline
 $H_{m,n}^{(Van)}$   &  1.18654642 &-0.23239211 & 0.08296761 & 0.13689074 \\
    \hline
  $H_{m,n}^{(Cor)}$  &  1.03589137 &-0.38304716 & 0.07197488 & 0.21862803 \\
\hline
$H_{m,n}^{(WG)}$   &  1.28060252 &-0.13833601 & 0.07087153  &0.08993751 \\
\hline
 $\widehat{H}_{\varepsilon;n}(X)$ & 1.33450763 &-0.08443091&  0.07002134 & \textbf{0.07707990}\\
    \hline
  \end{tabular}
  \caption{Results for  $n=10$, $m=3$, $h_n=0.157$, Normal distribution $N(0,1)$}\label{tab1}
\end{table}
\begin{table}[h!]
  \centering
  \begin{tabular}{|c|c|c|c|c|}
    \hline
     & Estimate & bias & Variance & MSE \\
    \hline
 $H_{m,n}^{(V)}$    &  1.10631536& -0.31262317 & 0.03182798 & 0.12952939\\
    \hline
 $H_{m,n}^{(Van)}$   &  1.24420922 &-0.17472931&  0.03532923 & 0.06582424 \\
    \hline
  $H_{m,n}^{(Cor)}$  &  1.24195337 &-0.17698517  &0.03230411 & 0.06359555\\
\hline
$H_{m,n}^{(WG)}$   &  1.36008222& -0.05885632 & 0.03182798  &0.03526021\\
\hline
 $\widehat{H}_{\varepsilon;n}(X)$ & 1.43214893& 0.01321040 &0.02920368 &\textbf{0.02934899}\\
    \hline
  \end{tabular}
  \caption{Results for  $n=20$, $m=3$, $h_n=0.081$, Normal distribution $N(0,1)$}\label{tab2}
\end{table}
\begin{table}[h!]
  \centering
  \begin{tabular}{|c|c|c|c|c|}
    \hline
     & Estimate & bias & Variance & MSE \\
    \hline
 $H_{m,n}^{(V)}$    & 1.26025240 &-0.15868613  &0.01126393 & 0.03643395\\
    \hline
 $H_{m,n}^{(Van)}$   &  1.29349009 &-0.12544844  &0.01210094  &0.02782615\\
    \hline
  $H_{m,n}^{(Cor)}$  &  1.35839575 &-0.06054278  &0.01148899  &0.01514293\\
\hline
$H_{m,n}^{(WG)}$   &  1.40287628 &-0.01606226  &0.01126393  &0.01151066\\
\hline
 $\widehat{H}_{\varepsilon;n}(X)$ & 1.42053167 &0.00159314 &0.01085020 &\textbf{0.01084189}\\
    \hline
  \end{tabular}
  \caption{Results for  $n=50$, $m=4$, $h_n=0.0333$, Normal distribution $N(0,1)$}\label{tab3}
\end{table}

\begin{table}[h!]
  \centering
  \begin{tabular}{|c|c|c|c|c|}
    \hline
     & Estimate & bias & Variance & MSE \\
    \hline
 $H_{m,n}^{(V)}$    & -0.142936202& -0.142936202 & 0.001730435 & 0.022161020\\
    \hline
 $H_{m,n}^{(Van)}$   & -0.000811817& -0.000811817 & 0.003428982 & 0.003429298\\
    \hline
  $H_{m,n}^{(Cor)}$  & -0.048455430& -0.048455430 & 0.001780981 & 0.004128732\\
\hline
$H_{m,n}^{(WG)}$   &-0.0003123278 &-0.0003123278 & 0.0017304350 & 0.0017303596\\
\hline
 $\widehat{H}_{\varepsilon;n}(X)$ & -0.001503714& -0.001503714&  0.001005806 & \textbf{0.001007967}\\
    \hline
  \end{tabular}
  \caption{Results for  $n=50$, $m=4$, $h_n=0.522$, Uniform distribution }\label{tab3U}
\end{table}

\begin{table}[h!]
  \centering
  \begin{tabular}{|c|c|c|c|c|}
    \hline
     & Estimate & bias & Variance & MSE \\
    \hline
 $H_{m,n}^{(V)}$    & -0.260890919 &-0.163204390 & 0.009863565 & 0.036497265\\
    \hline
 $H_{m,n}^{(Van)}$   & -0.20651205 &-0.10882552 & 0.01139025 & 0.02323097\\
    \hline
  $H_{m,n}^{(Cor)}$  &  -0.163758480&  -0.066071951 & 0.009953612 & 0.014317124\\
\hline
$H_{m,n}^{(WG)}$   & -0.118267044& -0.020580515 & 0.009863565 & \textbf{0.010285150}\\
\hline
 $\widehat{H}_{\varepsilon;n}(X)$ &-0.08547838&  0.01220815 & 0.01945246 & 0.01959761\\
    \hline
  \end{tabular}
  \caption{Results for  $n=50$, $m=4$, $h_n=0.6104$, Weibull distribution with
  the shape parameter  equal to $2$ and  the scale parameter equal to $0.5$}\label{tab3W}
\end{table}

\begin{table}[h!]
  \centering
  \begin{tabular}{|c|c|c|c|c|}
    \hline
     & Estimate & bias & Variance & MSE \\
    \hline
 $H_{m,n}^{(V)}$    &0.85568859 &-0.14431141 & 0.02233983  &0.04316114\\
    \hline
 $H_{m,n}^{(Van)}$   &  0.92199497 &-0.07800503 & 0.02313389 & 0.02921405\\
    \hline
  $H_{m,n}^{(Cor)}$  &  0.95566152 &-0.04433848  &0.02266589 & 0.02462726\\
\hline
$H_{m,n}^{(WG)}$   &  0.998312466 &-0.001687534&  0.022339826&  \textbf{0.022338205}\\
\hline
 $\widehat{H}_{\varepsilon;n}(X)$ & 1.31010002& 0.31010002 &0.06795533& 0.16410376\\
    \hline
  \end{tabular}
  \caption{Results for  $n=50$, $m=4$, $h_n=0.712$, Exponential distribution with parameter $1$}\label{tab3E}
\end{table}

\begin{table}[h!]
  \centering
  \begin{tabular}{|c|c|c|c|c|}
    \hline
     & Estimate & bias & Variance & MSE \\
    \hline
$H_{m,n}^{(V)}$    & 1.63721785 &-0.13625972&  0.02903662 & 0.04759752\\
    \hline
 $H_{m,n}^{(Van)}$   & 1.58175430 &-0.19172327 & 0.02360710 & 0.06036019\\
    \hline
  $H_{m,n}^{(Cor)}$  &  1.74658234 &-0.02689523 & 0.03048823  &0.03120548\\
\hline
$H_{m,n}^{(WG)}$   &  1.779841726& 0.006364154& 0.029036620 &0.029071315\\
\hline
 $\widehat{H}_{\varepsilon;n}(X)$ & 1.75882993 &-0.01464764 & 0.02592497 & \textbf{0.02613434}\\
    \hline
  \end{tabular}
  \caption{Results for  $n=50$, $m=4$, $h_n=0,0336$, Student's $t$-distribution
  with the number of degrees of freedom equal to $3$ }\label{tab3T3}
\end{table}

\begin{table}[h!]
  \centering
  \begin{tabular}{|c|c|c|c|c|}
    \hline
     & Estimate & bias & Variance & MSE \\
    \hline
 $H_{m,n}^{(V)}$    & 1.48096844 &-0.14653424  &0.02047265&  0.04194084\\
    \hline
 $H_{m,n}^{(Van)}$   &  1.46376243 &-0.16374024 & 0.01830929&  0.04511650\\
    \hline
  $H_{m,n}^{(Cor)}$  &1.58556888 & -0.04193379 & 0.02135329 & 0.02310746\\
\hline
$H_{m,n}^{(WG)}$   &  1.623592312 &-0.003910361 & 0.020472653&  0.020483849\\
\hline
 $\widehat{H}_{\varepsilon;n}(X)$ &1.621348438 &-0.006154234  &0.018522682&  \textbf{0.018556852}\\
    \hline
  \end{tabular}
  \caption{Results for  $n=50$, $m=4$, $h_n=0.344$, Student's $t$-distribution
   with the number of degrees of freedom equal to $5$ }\label{tab3T5}
\end{table}

\begin{table}[h!]
  \centering
  \begin{tabular}{|c|c|c|c|c|}
    \hline
     & Estimate & bias & Variance & MSE \\
    \hline
 $H_{m,n}^{(V)}$    & 2.51810441 &-0.01291983 & 0.09321374&  0.09336202\\
    \hline
 $H_{m,n}^{(Van)}$   &  2.24258072 &-0.28844353 & 0.05651327 & 0.13970163\\
    \hline
  $H_{m,n}^{(Cor)}$  &2.65247722 &0.12145298& 0.09936347 &0.11409442\\
\hline
$H_{m,n}^{(WG)}$   &  2.66072829 &0.12970404 &0.09321374 &0.11001824\\
\hline
 $\widehat{H}_{\varepsilon;n}(X)$ &2.49016605 &-0.04085819 & 0.07746897 &\textbf{ 0.07912286}\\
    \hline
  \end{tabular}
  \caption{Results for  $n=50$, $m=4$, $h_n=0.0235$,  Cauchy distribution}\label{tab3C5}
\end{table}

\noindent We turn now to compare robustness property of the above
estimators of the entropy $H(X)$.
The robustness here is to be stand versus contamination and not in
terms of influence function or break-down points which make sense
only under parametric or semiparametric settings.
According to \cite{Huber2009}: ``Most approaches to
 robustness are based on the following intuitive requirement: \emph{A discordant small 
 minority should never be able to override the evidence of the majority of the observations.}'' 
In the same reference, it is mentioned that resistant statistical procedure, i.e., if 
the value of the estimate is insensitive to small changes in the underlying sample, is 
equivalent to robustness for practical purposes in view of Hampel's theorem, refer to \cite[Section 1.2.3 and Section 2.6]{Huber2009}.
 Typical examples for the notion of robustness in the nonparametric setting are the 
sample mean and the sample median which are the nonparametric estimates of the population mean and median, respectively. Although nonparametric, the sample mean is highly    
sensitive to outliers and therefore for symmetric distribution and contaminated data the sample median is more appropriate to estimate the population mean or median, refer to \cite{Huber2009} for more details.\vskip5pt\noindent
In our simulation, we will consider data generated
from $N(0,1)$ distribution  where a ``small'' proportion
$\epsilon$ of observations were replaced by atypical ones
generated from a contaminating distribution $F^*(\cdot)$.  We
consider two cases with $\epsilon=4\%$ and $\epsilon=10\%$, and we
choose the contaminating distribution $F^*(\cdot)$ to be the
uniform distribution on $[0,1]$, the results are presented in
Table \ref{tab4} with $\epsilon=4\%$ and Table \ref{tab5} with $\epsilon=10\%$. Let $f(\cdot)$  denote the density function of $N(0,1)$ and $f^*(\cdot)$ the
density of the contaminating distribution $F^*(\cdot)$. The contaminated
sample, $X_1,\ldots,X_n,$ can be seen as if it has been
generated from the density $$f_\epsilon(\cdot):=(1-\epsilon)f(\cdot)+\epsilon
f^*(\cdot).$$ Since the sample is contaminated, all the above estimators
tend to $H(f_\epsilon)$ and not to $H(f)$.  The objective here is to
obtain the best (in the sense that the corresponding MSE is the
smallest) estimate of the entropy $H(f)$ from the contaminated
data $X_1,\ldots,X_n$. We will consider the five estimates as
above, and we compute their bias, variance and MSE by Monte-Carlo
using $5000$ replications. From Tables \ref{tab4}-\ref{tab5}, we can see  that our
estimator is the best one.

\begin{table}[h!]
  \centering
  \begin{tabular}{|c|c|c|c|c|}
    \hline
     & Estimate & bias & Variance & MSE \\
    \hline
 $H_{m,n}^{(V)}$    & 1.24299668 &-0.17594185  &0.01115852  &0.04210290\\
    \hline
 $H_{m,n}^{(Van)}$   &  1.27422798 &-0.14471055  &0.01143207  &0.03236179\\
    \hline
  $H_{m,n}^{(Cor)}$  &  1.34170579 &-0.07723274 & 0.01160997  &0.01756326\\
\hline
$H_{m,n}^{(WG)}$   &  1.38562055 &-0.03331798  &0.01115852  &0.01225745\\
\hline
 $\widehat{H}_{\varepsilon;n}(X)$ & 1.40066144 &-0.01827709  &0.01054644  &\textbf{0.01086995}\\
    \hline
  \end{tabular}
  \caption{Results for  $n=50$, $m=4$, $h_n=0.0333$, $\epsilon=4\%$, Normal distribution $N(0,1)$}\label{tab4}
\end{table}

\begin{table}[h!]
  \centering
  \begin{tabular}{|c|c|c|c|c|}
    \hline
     & Estimate & bias & Variance & MSE \\
    \hline
 $H_{m,n}^{(V)}$    & 1.22256006 &-0.19637848 & 0.01245586  &0.05100791\\
    \hline
 $H_{m,n}^{(Van)}$   &  1.24540126 &-0.17353727 & 0.01355086 & 0.04365249\\
    \hline
  $H_{m,n}^{(Cor)}$  &  1.32178766 &-0.09715087 & 0.01243944&  0.02186529\\
\hline
$H_{m,n}^{(WG)}$   &  1.36518393 &-0.05375460 & 0.01245586  &0.01533296\\
\hline
 $\widehat{H}_{\varepsilon;n}(X)$ & 1.37585099 &-0.04308754 & 0.01219767 & \textbf{0.01404201}\\
    \hline
  \end{tabular}
  \caption{Results for  $n=50$, $m=4$, $h_n=0.0333$, $\epsilon=10\%$, Normal distribution $N(0,1)$}\label{tab5}
\end{table}

\begin{remark}
To understand well the behavior of the proposed estimator, it will be interesting
to consider several family of distributions:
\begin{enumerate}
\item Distribution with support $(-\infty,\infty)$ and  symmetric.
\item Distribution with  support $(-\infty,\infty)$ and asymmetric.
\item Distribution with  support $(0,\infty)$.
\item Distribution with  support $]0,1[$.
\end{enumerate}
Extensive simulations for these families will be undertaken elsewhere.
\end{remark}
In the following remark, we give a way how the choose the smoothing parameter in practice.

\begin{remark}
The choice of the smoothing parameter plays an instrumental role in implementation of
practical estimation. We recall that the smoothing parameter $h_{n}$, in our simulations,
 has been chosen to minimize the MSE of the estimator, assuming that the underlying
  distribution is known. In more general case, without assuming any knowledge on the
  underlying distribution function, one can use, among others, the selection procedure 
  proposed in \cite{Jones1992}. \cite{Jones1992} derived that the asymptotic MSE
  of $\widehat{q}_{n}(\cdot)$, in the case of convolution-kernel estimator, which is given by
\begin{eqnarray*}
AMSE(\widehat{q}_{n}(t))&=&\frac{h_{n}^{4}}{4}{q^{\prime\prime}}^{2}(t)\left
\{\int_{\mathbb{R}}x^{2}K(x)dx\right\}^{2}\\&& +\frac{1}{nh_{n}}q^2(t)\int_\mathbb{R}K^{2}(x) dx.
\end{eqnarray*}
Minimizing the last equation with respect to $h_{n}$, we find the asymptotically
 optimal bandwidth for $\widehat{q}_{n}(\cdot)$ as
\begin{equation*}
h_{n}^{\rm opt}=\left\{\frac{q(t)^{2}\int_\mathbb{R} K^{2}(x)dx}{n({q^{\prime\prime}}
^{2}(t))\left\{\int_{\mathbb{R}}x^{2}K(x)dx\right\}^{2}}\right\}^{1/5}.
\end{equation*}
Note that $ h_{n}^{\rm opt} $ depends on the unknown functions
$$
q(t)=\frac{1}{f(Q(t))},
$$
and
$$
q^{\prime\prime}(t)=\frac{f^{\prime\prime}(Q(t))f(Q(t))-3[f^{\prime}(Q(t))]^{2}}{f^{5}(Q(t))}.
$$
These functions may be estimated in the classical way,  refer to \cite{Jones1992},
 \cite{Silverman1986} and \cite{Cheng2006} for more details on the subject and the
 references therein. Another way to estimate the optimal value of $h_{n}$ is to
 use a cross-validation type method.
\end{remark}

\begin{remark}
 The main problem in using entropy estimates
such as in (\ref{6-entropy-estimator}) is to choose properly the smoothing
parameter $h_n$.  With a lot more effort, we could derive analog results here
 for $\widehat{H}_{\varepsilon;n}(X)$  using
the methods in \cite{Bouzebda-Elhattab2009,Bouzebda-Elhattab2010,Bouzebda2011},
as well as the modern empirical process tools developed
in \cite{Mason2005} in their work on uniform
in bandwidth consistency of kernel-type estimators.
\end{remark}

\section{Test for normality}\label{simulation}
\noindent A well-known theorem of information theory [see, e.g.,
p. 55, \cite{shannon1948}] states that among all distributions
that possess a density function $f(\cdot)$ and have a given
variance $\sigma^2$, the entropy $H(X)$ is maximized by the normal
distribution. The entropy of the normal distribution with variance
$\sigma^2$ is $\log\big(\sigma\sqrt{2\pi e}\big)$. As pointed out
by \cite{vasicek1976}, this property can be used for tests of
normality. Towards this aim, one can use the estimate $\widehat{H}_{\varepsilon;n}(X)$ of
$H(X)$, as follows. Let $X_1,\ldots,X_n$ be independent random
replic{\ae} of a random variable $X$ with quantile density
function $q(\cdot)$. Let $\widehat{H}_{\varepsilon;n}(X)$ the estimator of $H(X)$ as in
(\ref{6-entropy-estimator}). Let $\mathds{T}_n$ denote the statistic

\begin{equation}\label{ourstat}
\mathds{T}_n:=\log\left(\sqrt{2\pi \sigma_n^{2}}\right)+0.5-\widehat{H}_{\varepsilon;n}(X),
\end{equation}
where
$\sigma_n^2$ is the sample standard deviation based on
$X_1,\ldots,X_n,$ defined by
$$\sigma_n^2:=\frac{1}{n}\sum_{i=1}^{n}\left(X_i-\overline{X}_n\right)^2,$$
 where $$\overline{X}_n:=n^{-1}\sum_{i=1}^{n}X_i.$$ The normality
hypothesis will be rejected whenever the observed value of $\mathds{T}_n$
will be significantly less than $0$,
in the sense we precise below. The exact critical values
$T_\alpha$ of $\mathds{T}_n$ at the significance levels
$\alpha\in]0,1[$ are defined by the equation
\begin{equation*}
    \mathbb{P}\left(\mathds{T}_n\leq T_\alpha\right)=\alpha.
\end{equation*}
The distribution of $\mathds{T}_n$ under the null hypothesis
cannot readily be obtained analytically. To evaluate the critical
values $T_\alpha$, we have used a Monte Carlo method, for sample
sizes $10\leq n\leq 50$ and the significance value given by
$\alpha=0.05$. Namely, for each $n\leq 50$, we simulate 20000
samples of size $n$ from the standard normal distribution. Since
$\alpha=0.05=1000/20000$, we determine the 1000-th order statistic
$t_{1000,20000}$ and obtain the critical value $T_{n,0.05}$ through
the equation $T_{n,0.05}=t_{1000,20000}$.
Our results are presented in Table \ref{tableau}.\\
\begin{table}[h!]
  \centering
\begin{tabular}{cccccc}
\hline\hline
\multicolumn{1}{c}{Sample size} &\multicolumn{5}{c}{ Percentage level } \\
\cline{2-6}
\multicolumn{1}{c}{$n$ } &  0.1 & 0.05 & 0.025 & 0.01 & 0.005 \\
\hline \hline
35 & 0.03660258& 0.05896114 &0.07819581&0.1041396 &0.1229790\\
40 &  0.03641781 & 0.05732028&0.07655416& 0.1003001&0.1177802 \\
45 &  0.03011983 &0.04910612 & 0.06554724& 0.0844859& 0.1004404 \\
50 &  0.02534047 &0.04232442 &0.05874272 &0.07830533 &0.09371921 \\
 \hline
\end{tabular}
\caption{Critical Values of $ \mathds{T}_{n}$.}\label{tableau}
\end{table}\vskip7pt

\noindent\cite{Park2003} establish the entropy-based goodness of fit test statistics based on
the nonparametric distribution functions of the sample entropy and modified sample entropy
\cite{Ebrahimi1994}, and compare their performances for the exponential and normal distributions.
 The authors consider
$$
H_{m,n}^{Ebr}:=\frac{1}{n}\sum_{i=1}^{n}\log \left(\frac{n}{c_{i}m}\{X_{i+m;n}-X_{i-m;n}\}\right),
$$
where
$$c_{i}:=\left\{
\begin{array}{lcl}
     1+\frac{i-1}{m} &\mbox{ if } & 1\leq i\leq m ,  \\
      2 &\mbox{ if } & m+1\leq i\leq n-m,   \\
     1+\frac{n-i}{n} & \mbox{ if } & n-m\leq i \leq n.
\end{array}\right.
$$
 \cite{Yousefzadeh2008} use a new cdf estimator to obtain a nonparametric entropy estimate
 and use it for testing exponentially
  and normality, they introduce the following estimator
$$
H_{m,n}^{You}:=\sum_{i=1}^{n}\log \left(\frac{X_{i+m;n}-X_{i-m;n}}{\widehat{F}_{n}(X_{i+m;n})-
\widehat{F}_{n}(X_{i-m;n})}\right)\frac{\widehat{F}_{n}(X_{i+m;n})-\widehat{F}_{n}(X_{i-m;n})}
{\sum_{i=1}^{n}(\widehat{F}_{n}(X_{i+m;n})-\widehat{F}_{n}(X_{i-m;n}))},
$$
 where

 $$\widehat{F}_{n}(x):=\left\{
 \begin{array}{ lcl}
    \frac{n-1}{n(n+1)}\left(\frac{n}{n-1}+\frac{x-X_{0;n}}{X_{2;n}-X_{0;n}}+\frac{x-X_{1;n}}{X_{3;n}-X_{1;n}}\right)
     &\mbox{for}&     X_{1;n}\leq x\leq X_{2;n},\\
    &&\\
       \frac{n-1}{n(n+1)}\left(i+\frac{1}{n-1}+\frac{x-X_{i-1;n}}{X_{i+1;n}-X_{i-1;n}}+\frac{x-X_{i;n}}
       {X_{i+2;n}-X_{i;n}}\right)&\mbox{for}& X_{i;n}\leq x\leq X_{i+1;n},  \\
       && i=2,\ldots,n-2,\\
       \frac{n-1}{n(n+1)}\left(n-1\frac{1}{n-1}+\frac{x-X_{n-2;n}}{X_{n;n}-X_{n-2;n}}+\frac{x-X_{n-1;n}}
       {X_{n+1;n}-X_{n-1;n}}\right) &\mbox{for}&
          X_{n-1;n}\leq x\leq X_{n,n},
\end{array}\right.
 $$
 and
 \begin{eqnarray*}
 X_{0;n}&:=&X_{1;n}-\frac{n}{n-1}(X_{2;n}-X_{1;n}),\\
 X_{n+1;n}&:=&X_{n;n}-\frac{n}{n-1}(X_{n;n}-X_{n-1;n}).
 \end{eqnarray*}
To compare the power of the proposed test a simulation was conducted. We used tests based on the following entropy
estimators: $H_{m,n}^{(WG)}$,  $H_{m,n}^{Ebr}$, $H_{m,n}^{You}$ and our test  (\ref{ourstat}).  In the power comparison, we consider the following alternatives.
\begin{enumerate}
\item[$\bullet$] The Weibull distribution with density function
$$
f(x;\lambda, k)=\frac{k}{\lambda}\left(\frac{x}{\lambda}\right)\exp\left(\left(\frac{-x}{\lambda}\right)^{k}\right)\mathds{1}\{x>0\},
$$
where $k > 0$ is the shape parameter and $\lambda>0$ is the scale parameter of the distribution
and where $\mathds{1}\{\cdot\}$ stands for the indicator function of the set $\{\cdot\}$.
\item[$\bullet$] The uniform distribution with density function
$$
f(x)=1,~~0<x<1.
$$
\item[$\bullet$] The Student $t$-distribution with density function
$$
f(x;\nu)=\frac{\Gamma((\nu+1)/2)}{\Gamma(\nu/2)}\frac{1}{\sqrt{\nu \pi}}\frac{1}{(1+x^{2}/\nu)^{(\nu+1)/2}},~~\nu>2,~~-\infty<x<\infty.
$$
\end{enumerate}

\begin{table}[h!]
  \centering
\begin{tabular}{cccccc}
\hline\hline
&\multicolumn{5}{c}{Statistics based on } \\
\cline{2-6}
\multicolumn{1}{c}{Alternatives } & $\widehat{H}_{\varepsilon;n}(X)$ & $H_{4,n}^{(WG)}$ &$H_{8,n}^{(WG)}$& $H_{m,n}^{Ebr}$ &$H_{m,n}^{You}$ \\
\hline \hline
Uniform & 0.9999 & 0.9262 & 0.96685 &  0.9275  & 0.8768 \\
Weibull(2) & 0.8264& 0.3297&0.33795 &  0.4211 & 0.3444 \\
Student $t_{5}$ & 0.9306 & 0.1358 &0.05530 & 0.1484 & 0.2345\\
Student $t_{3}$ & 1.0000 & 0.3696 & 0.18245& 0.3736&0.5124 \\
\hline
\end{tabular}
\caption{Power estimate of $0.05$ tests against alternatives of the normal distribution
based on 20 000 replications for sample size $n$=50.}\label{tableau23}
\end{table}\vskip7pt
\begin{table}[h!]
  \centering
\begin{tabular}{ccc}
\hline\hline
\multicolumn{1}{c}{Alternatives } &$h_{n}$  \\
\hline \hline
Uniform & $h_{n}=0.5297$ \\
Weibull(2) &  $h_{n}=0.6555$   \\
Student $t_{5}$ &$h_{n}=0.0310$  \\
Student $t_{3}$ & $h_{n}=0.0189$   \\
\hline
\end{tabular}
\caption{Choice of smoothing parameter for $\widehat{H}_{\varepsilon;n}(X)$.}\label{tableau235}
\end{table}\vskip7pt

\begin{remark}
We mention the value taken in Table \ref{tableau23} for the statistics based on $H_{m,n}^{Ebr}$ and $H_{m,n}^{You}$
are the same to that calculated in Table 6, p. 1493 of \cite{Yousefzadeh2008}, for the same alternatives that we consider in our comparison.
\end{remark}
In Table \ref{tableau23} we have reported the results of power comparison for the  sample size is 50.
We made $20000$ Monte-Carlo simulations to compare the powers of the proposed test against $4$
 alternatives.
From Table \ref{tableau23}, we can see that the proposed test $\mathds{T}_{n}$ shows
better power than all other statistics for the alternative that we consider.  It is
 natural that our test has good performances for unbounded support since the kernel-type estimators behave well in this situation.

\section{Concluding remarks and future works}\label{conclusion}
\noindent We have proposed a new estimator of entropy based on the
kernel-quantile density estimator. Simulations show that this
estimator behaves better than the other competitors both under
contamination or not. More precisely, the MSE is consistently
smaller than the MSE of the spacing-based estimators considered in
our study. It will be interesting to compare theoretically the power and the
robustness of the test of normality, based on the proposed estimator of
the present paper,  with those considered
in \cite{Esteban2001}. The study of entropy in presence of censored data
is possible using a similar methodology as presented here. It would be interesting
to provide a complete investigation of the choice of the parameter $h_{n}$ for kernel difference-type estimator which requires nontrivial mathematics, this
would go well beyond the scope of the present paper. \vskip5pt\noindent
The problems and the methods described here all are inherently univariate. A natural and useful multivariate extension
is the use of copula function. We propose to extend the results of this paper
to the multivariate case in the following way. Consider a $\mathbb{R}^d$-valued random vector $\mathbf{X }= (X_1,\ldots,X_d)$ with
joint cdf $$\mathbb{F}(\mathbf{x}):=\mathbb{F}(x_{1},\ldots,x_{d}):=\mathbb{P}(X_{1}
\leq x_{1},\ldots,X_{d}\leq x_{d})$$ and
marginal cdf's $$F_j(x_{j}):=\mathbb{P}(X_{j}\leq x_{j})~~\mbox{ for }~~j=1,\ldots,d.$$ If the marginal distribution functions
$F_1(\cdot),\ldots, F_d(\cdot)$ are continuous,
then, according to Sklar's theorem [\cite{Sklar1959}], the copula function $\mathbb{C}(\cdot)$, pertaining to $\mathbb{F}(\cdot)$, is unique and
$$\mathbb{C}(\mathbf{u}):=\mathbb{C}(u_1,\ldots, u_d) := \mathbb{F}(Q_1(u_1),\ldots, Q_d (u_d)),  ~~ \mbox{ for } ~~\mathbf{u}\in [0,1]^{d},$$
where, for $j = 1,\ldots, d$, $Q_j (u) := \inf\{x : F_j(x)\geq u\}$  with $u\in(0, 1]$, $Q_j (0):=\lim_{t\downarrow 0}Q_{j}(t):=Q_j (0^{+})$, is
 the quantile function of $F_j(\cdot)$.
The differential entropy may represented via density  copula $ \mathbf{c}(\cdot)$ and
quantile densities $q_i(\cdot)$  as follows
\begin{equation}\label{repentropy}
H(\mathbf{X}) = \sum_{i=1}^d H(X_i)+ H(\mathbf{c}),
\end{equation}
where
\begin{equation}\label{entropyuniv}
H(X_i) = \int_{[0,1]}\log\big(q_{i}(u)\big) du,
\end{equation}
 $q_{i}(u)=dQ_{i}(u)/du$, for $i=1,\ldots,d$
and $H(\mathbf{c})$ is the \emph{copula entropy} defined by
\begin{eqnarray}\label{entropy}
\nonumber H(\mathbf{c}) &=& -\int_{\mathbb{R}^d} \mathbf{c}(F_1(x_1),\ldots,F_d(x_d))\log\big(\mathbf{c}(F_1(x_1),\ldots,F_d(x_d))\big)
d{\bf x}\\\nonumber&=&-\int_{[0,1]^d} \mathbf{c}(u_1,\ldots,u_d)\log\big(\mathbf{c}(u_1,\ldots,u_d)\big)
d{\bf u}\\&=&-\int_{[0,1]^d} \mathbf{c}(\mathbf{u})\log\big(\mathbf{c}(\mathbf{u})\big)
d{\bf u},
\end{eqnarray}
where $$\mathbf{c}(u_1,\ldots,u_d) = \frac{\partial^{d}}{\partial u_1\ldots\partial u_d
}\mathbb{C}(u_1,\ldots,u_d)$$ is the copula density.

\section{Proof.}\label{6-demo}
This section is devoted to the proofs of our results.  The previously defined notation continues to be used below.\\
\noindent {\bf Proof of Theorem \ref{6-theorem1}.}
Recall that  we set for all $\varepsilon\in]0,1/2[$, $U(\varepsilon)=[\varepsilon,1-\varepsilon]$, and
\begin{equation*}
    H_\varepsilon(X)=\varepsilon\log\big(q(\varepsilon)\big)+
    \int_{U(\varepsilon)}\log\big(q(x)\big)dx+\varepsilon\log\big(q(1-\varepsilon)\big).
\end{equation*}
  Therefore,  by (\ref{Keyidea}),  we have the following
\begin{eqnarray}\label{6-decompositio}
\nonumber \left|\widehat{H}_{\varepsilon;n}(X)-H(X)\right| &\leq&
\left| \widehat{H}_{\varepsilon;n}(X)-H_\varepsilon(X)\right|+\left|H_\varepsilon(X)-H(X)\right|
\\&=&\left| \widehat{H}_{\varepsilon;n}(X)-H_\varepsilon(X)\right|+o\big(\eta(\varepsilon)\big).
\end{eqnarray}
We evaluate the first term of the right hand of (\ref{6-decompositio}). By the triangular inequality, we have
\begin{eqnarray}
\nonumber\left| \widehat{H}_{\varepsilon;n}(X)-H_\varepsilon(X)\right|&=&\left|
\varepsilon\left(\log\big(\widehat{q}_n(\varepsilon)\big)-\log\big(q(\varepsilon)\big)\right)+\int_{U(\varepsilon)}
(\log\big(\widehat{q}_n(x)\big)-\log\big(q(x)\big))dx\right.
\\&&\nonumber+\varepsilon\left(\log\big(\widehat{q}_n(1-\varepsilon)\big)-\log\big(q(1-\varepsilon)\big)\right)
\Big|\\&\leq&\nonumber
\left|
\varepsilon\left(\log\big(\widehat{q}_n(\varepsilon)\big)-\log\big(q(\varepsilon)\big)\right)\right|\\&&\nonumber+\left|\int_{U(\varepsilon)}
(\log\big(\widehat{q}_n(x)\big)-\log\big(q(x)\big))dx\right|
\\&&\nonumber+\left|\varepsilon\left(\log\big(\widehat{q}_n(1-\varepsilon)\big)-\log\big(q(1-\varepsilon)\big)\right)
\right|.
\end{eqnarray}
 We note that for
$z>0$, $$|\log z|\leq |z-1|+|1/z-1|.$$ Therefore, we have
\begin{eqnarray*}
\lefteqn{\left|\log\left(\widehat{q}_n(x)\right)-\log\big(q(x)\big)\right|
=\left|\log\left(\frac{\widehat{q}_n(x)}{q(x)}\right)\right|}\\
&\leq& \frac{\left|\widehat{q}_n(x)-q(x)\right|}{q(x)}
+\frac{\left|\widehat{q}_n(x)-q(x)\right|}{\widehat{q}_n(x)}.
\end{eqnarray*}
Under conditions (Q.1-2-3) and (K.1-2-3), we may apply Theorem 2.2 in
\cite{chengParzen1997}, for all fixed $\varepsilon \in ]0,1/2[$ and recall $M(q;K_n)$ defined in Theorem \ref{6-theorem1},
 we have, as $n\rightarrow\infty$,
\begin{equation}\label{chenpar}
\sup_{x\in
U(\varepsilon)}\left|\widehat{q}_n(x)-q(x)\right|=O_{\mathbb{P}}\Big(n^{-1/2}M(q;K_n)+n^{-\beta}\Big),
\end{equation}
we infer that we have, uniformly over
$x\in U(\varepsilon)$, $\widehat{q}_n(x)\geq (1/2)q(x),$ for all $n$ enough
large. This fact implies
\begin{eqnarray}
\nonumber\left| \widehat{H}_{\varepsilon;n}(X)-H_\varepsilon(X)\right|&\leq&4\sup_{x\in U(\varepsilon)}\left|\widehat{q}_n(x)-q(x)\right|,
\end{eqnarray}
in probability,
which implies that
\begin{eqnarray}
\nonumber\left| \widehat{H}_{\varepsilon;n}(X)-H_\varepsilon(X)\right|=O_{\mathbb{P}}\Big(n^{-1/2}M(q;K_n)+n^{-\beta}+\eta(\varepsilon)\Big).
\end{eqnarray}
Thus the proof is complete.\hfill$\Box$\\
\noindent {\bf Proof of Theorem \ref{6-norm1}.}  Throughout, we will make use
of the fact (see e.g.  \cite{CP1978,csorgorevesz1981}) that we can
define $X_1,\ldots,X_n$ on a probability space which carries a
sequence of Brownian bridges $$\{B_n(t):t\in[0,1], n\geq 1\},$$ such that,
\begin{eqnarray}\label{qqqq2}
\sqrt{n}(Q_n(t)-Q(t))-q(t)B_n(t)=q(t)e_n(t)
\end{eqnarray}
and
\begin{equation}
\lim\sup_{n\rightarrow \infty}\left[\frac{n^{1/2}}{\mathbb{A}_{\gamma}(n)}\right]q(t)e_n(t)\leq \mathbb{C},
\end{equation}
with probability one, where $\mathbb{C}$ is a universal constant and
\begin{eqnarray*}
\mathbb{A}_{\gamma}(n):=
\left\{
\begin{array}{ll}
 \log n,&\max\{q(0),q(1)\}<\infty ~~\mbox{or}~~\gamma\leq 2,\\
 (\log\log n)^{\gamma}(\log n)^{(1+\nu)(1-\gamma)}& \gamma>2,
\end{array}\right.
\end{eqnarray*}
with $\gamma$ in (Q.2) and an arbitrary positive $\nu$. Using Taylor expansion we get, for  $\lambda \in ]0,1[$,
\begin{eqnarray}\label{6-tilds}
\nonumber \lefteqn{\sqrt{n}(\widehat{H}_{\varepsilon;n}(X) -H_\varepsilon(X))}\\\nonumber& = &
\int_{U(\varepsilon)}\sqrt{n}(\log\big(\widehat{q}_n(x)\big)-\log\big(q(x)\big)) dx\\
&&\nonumber+\sqrt{n}\varepsilon(\log\big(\widehat{q}_n(\varepsilon)\big)-\log\big(q(\varepsilon\big))\\&&\nonumber
 +\sqrt{n}\varepsilon(\log\big(\widehat{q}_n(1-\varepsilon)\big)-\log\big(q(1-\varepsilon)\big)) \\
\nonumber  & = &  \int_{U(\varepsilon)}\sqrt{n}\left(\frac{\widehat{q}_n(x)-q(x)}{\lambda\widehat{q}_n(x)
+(1-\lambda)q(x)}\right) dx+\mathbb{L}_{n}(\varepsilon)\\
&=:& \mathbb{T}_{n}(\varepsilon)+\mathbb{L}_{n}(\varepsilon).
\end{eqnarray}
We have, under our assumptions, $\widehat{q}_n(x)\rightarrow q(x)$ in probability, which implies that
$\lambda\widehat{q}_n(x)+(1-\lambda)q(x)\rightarrow q(x)$ in probability.
Thus, we have, by Slutsky's theorem, as $n$ tends to the infinity, $\mathbb{T}_{n}(\varepsilon)$
 has the same limiting law of
\begin{equation}
\widetilde{\mathbb{T}_{n}}(\varepsilon):=\int_{U(\varepsilon)}(1/q(x))\sqrt{n}\left(\widehat{q}_n(x)-q(x)\right) dx.
\end{equation}  We have the following decomposition
\begin{eqnarray*}
\widetilde{\mathbb{T}_{n}}(\varepsilon)& = & \nonumber\int_{U(\varepsilon)}(1/q(x))\sqrt{n}\frac{d}{dx}
\left(\int_0^1Q_n(v)K_n\left(x,v\right)d\mu_n(v)-Q(x)\right) dx
\\& = & \nonumber\int_{U(\varepsilon)}(1/q(x))\frac{d}{dx}\left(\int_0^1\sqrt{n}(Q_n(v)-Q(v))K_n
\left(x,v\right)d\mu_n(v)\right) dx\\&  &- \nonumber\int_{U(\varepsilon)}(1/q(x))\sqrt{n}\frac{d}{dx}
\left(\int_0^1Q(v)K_n\left(x,v\right)d\mu_n(v)-Q(x)\right) dx
\end{eqnarray*}
Using   condition (K.3) and $$\int_{U(\varepsilon)}(1/q(x))dx\leq 1,$$ we have the following
\begin{eqnarray}
\left|\nonumber\int_{U(\varepsilon)}(1/q(x))\sqrt{n}\frac{d}{dx}\left(\int_0^1Q(v)K_n\left(x,v\right)d\mu_n(v)-Q(x)\right) dx\right|=O(n^{1/2-\beta}).
\end{eqnarray}
This, in turn, implies, in connection with (\ref{qqqq2}), that
\begin{eqnarray*}
\widetilde{\mathbb{T}_{n}}(\varepsilon)&=& \nonumber\int_{U(\varepsilon)}(1/q(x))\frac{d}{dx}
\left(\int_0^1(q(v)B_{n}(v)+ q(v)e_{n}(v))K_n\left(x,v\right)d\mu_n(v)\right) dx\\&&+O(n^{1/2-\beta}).
\end{eqnarray*}
Recalling the arguments of the proof of Lemma 3.2. of \cite{cheng1995}, 
for $$I_{n}(u)=[u-\delta_{n},u+\delta_{n}] ~~\mbox{and}~~ I_{n}^{c}(u)=[0,1]\backslash I_{n}(u),$$
we can show that
\begin{eqnarray*}
\lefteqn{\left|\nonumber\int_{U(\varepsilon)}(1/q(x))\frac{d}{dx}\left(\int_0^1q(v)e_{n}
(v)K_n\left(x,v\right)d\mu_n(v)\right) dx\right|}\\&\leq&\nonumber\int_{U(\varepsilon)}
(1/q(x))dx\sup_{x\in U(\varepsilon)}\left|\frac{d}{dx}\left(\int_0^1q(v)e_{n}(v)K_n\left(x,v\right)d\mu_n(v)\right) \right|.
\end{eqnarray*}
Then, it follows
\begin{eqnarray*}
\lefteqn{\left|\nonumber\int_{U(\varepsilon)}(1/q(x))\frac{d}{dx}\left(\int_0^1q(v)e_{n}(v)K_n\left(x,v\right)d\mu_n(v)\right) dx\right|}\\
&\leq&\sup_{x\in U(\varepsilon)}\left(\int_0^1q(v)\left|e_{n}(v)\frac{d}{dx}K_n\left(x,v\right)\right|d\mu_n(v)\right)
\\&\leq&\sup_{x\in U(\varepsilon)}|e_{n}(v)|\left(\int_0^1q(v)\left|\frac{d}{dx}K_n\left(x,v\right)\right|d\mu_n(v)\right)
\\&\leq&\mathbb{C}n^{-1/2}\mathbb{A_{\gamma}}(n)\left(\int_0^1q(v)\left|\frac{d}{dx}K_n\left(x,v\right)\right|d\mu_n(v)\right)
\\&\leq&\mathbb{C}n^{-1/2}\mathbb{A_{\gamma}}(n)\sup_{x\in U(\varepsilon-\delta_{n})}|q(x)|\sup_{x\in U(\varepsilon)}\left(\int_{I_{n}(x)}\left|\frac{d}{dx}K_n\left(x,v\right)\right|d\mu_n(v)\right) \\&&+
\mathbb{C}n^{-1/2}\mathbb{A_{\gamma}}(n)\sup_{x\in U(\varepsilon)}\int_{I_{n}^{c}(x)}\left|\frac{d}{dx}q(v)K_n\left(x,v\right)\right|d\mu_n(v)\\&=:&
\mathbb{R}_{n}^{(1)}(\varepsilon)+\mathbb{R}_{n}^{(2)}(\varepsilon).
\end{eqnarray*}
Since $q(\cdot)$
 is continuous on $]0,1[$, then it is bounded on compact interval $U(\varepsilon-\delta_{n})\subset ]0,1[$, which gives, by condition (K.2)
\begin{eqnarray*}
\mathbb{R}_{n}^{(1)}(\varepsilon)&\leq& \mathbb{C}n^{-1/2}\mathbb{A_{\gamma}}(n)\sup_{x\in U(\varepsilon-\delta_{n})}|q(x)|\\&=&
\mathbb{C}n^{-1/2}\mathbb{A_{\gamma}}(n)O(1)=o(1)\\&=&o_{\mathbb{P}}(1).
\end{eqnarray*}
Let $s=1+\delta$ with $\delta>0$ arbitrary, and let $r=1-1/s$. Then H\"older's inequality implies
\begin{eqnarray*}
\mathbb{R}_{n}^{(2)}(\varepsilon)&\leq& \mathbb{C}n^{-1/2}\mathbb{A_{\gamma}}(n)\sup_{x\in U(\varepsilon)}
\left\{\left[\int_{[0,1]}[q(v)]^{r}\frac{d}{dx}K_n\left(x,v\right)d\mu_n(v)\right]^{1/r}
\right.\\&&\times
\left.\left[\int_{[0,1]}\mathds{1}^{s}_{I_{n}^{c}(x)}(v)\frac{d}{dx}K_n\left(x,v\right)d\mu_n(v)\right]^{1/s}\right\}.
\end{eqnarray*}
Under condition (K.3), we can see that
\begin{eqnarray*}
\lefteqn{\sup_{x\in U(\varepsilon)}
\left[\int_{[0,1]}[q(v)]^{r}\frac{d}{dx}K_n\left(x,v\right)d\mu_n(v)\right]^{1/r}}\\
&\leq& \sup_{x\in U(\varepsilon)}\left|q(x)-
\left[\int_{[0,1]}[q(v)]^{r}\frac{d}{dx}K_n\left(x,v\right)d\mu_n(v)\right]^{1/r}\right|\\&=&O(n^{-\beta}).
\end{eqnarray*}
This in connection with condition (K.2), implies
\begin{eqnarray*}
\mathbb{R}_{n}^{(2)}(\varepsilon)&=&\mathbb{C}n^{-1/2}\mathbb{A_{\gamma}}(n)O(1)O\left(\sup_{x\in
 U(\varepsilon)}R_{n}^{1/(1/\delta)}\right)=o(1)\\&=& o_{\mathbb{P}}(1).
\end{eqnarray*}
Hence
\begin{eqnarray*}
\left|\nonumber\int_{U(\varepsilon)}(1/q(x))\frac{d}{dx}\left(\int_0^1q(v)e_{n}(v)K_n
\left(x,v\right)d\mu_n(v)\right) dx\right|=o_{\mathbb{P}}(1).
\end{eqnarray*}
Then, we conclude that
\begin{eqnarray*}
\widetilde{\mathbb{T}_{n}}(\varepsilon)& = & \nonumber\int_{U(\varepsilon)}(1/q(x))\frac{d}{dx}
\left(\int_0^1q(v)B_{n}(v)K_n\left(x,v\right)d\mu_n(v)\right) dx\\&&+o_{\mathbb{P}}(1).
\end{eqnarray*}
We have by \cite{chengParzen1997}, p. 297, (we can refer also to \cite{Stadtm1,Stadtm2} and \cite{Xiang1994}
for related results on smoothed Wiener process and Brownian bridge)
\begin{eqnarray*}
\lefteqn{\sup_{x\in U(\varepsilon)}\left|\int_0^1q(v)B_{n}(v)K_n\left(x,v\right)d\mu_n(v) -q(x)B_{n}(x)\right|}\\&\leq&
\sup_{x\in U(\varepsilon)}\int_0^1q(v)\left|B_{n}(v)-B_{n}(x)\right|K_n\left(x,v\right)d\mu_n(v)\\&&+
\sup_{x\in U(\varepsilon)}\left|B_{n}(x)\int_0^1\left|q(v)-q(x)\right|K_n\left(x,v\right)d\mu_n(v)\right|\\
&=&O_{\mathbb{P}}((2\delta_{n}\log \delta_{n}^{-1})^{1/2}))+O_{\mathbb{P}}\left(R_{n}^{1/(1+\delta)}\right)
+O_{\mathbb{P}}(n^{-\beta})
\\&=&o_{\mathbb{P}}(1).
\end{eqnarray*}
This gives, under condition (Q.1),
\begin{eqnarray*}
 \lefteqn{\nonumber\int_{U(\varepsilon)}(1/q(x))\frac{d}{dx}\left(\int_0^1q(v)B_{n}(v)K_n\left(x,v\right)d
 \mu_n(v)\right) dx}\\&=&(1/q(x))\left(\int_0^1q(v)B_{n}(v)K_n\left(x,v\right)d\mu_n(v)\right)
 \Big|_{\varepsilon}^{1-\varepsilon}\\&&-\int_{U(\varepsilon)}
 \left(\frac{d}{dx}(1/q(x))\right)\left(\int_0^1q(v)B_{n}(v)K_n\left(x,v\right)d\mu_n(v)\right) dx
 \\&=&(1/q(x))q(x)B_{n}(x) \Big|_{\varepsilon}^{1-\varepsilon}-\int_{U(\varepsilon)}\left(\frac{d}
 {dx}(1/q(x))\right)q(x)B_{n}(x) dx+o_{\mathbb{P}}(1)\\&=&\int_{U(\varepsilon)}(q(x)^{'}/q(x))B_{n}(x)
  dx+B_{n}(1-\varepsilon)-B_{n}(\varepsilon)+o_{\mathbb{P}}(1).
\end{eqnarray*}
Making use of \cite[Theorem 1.4.1]{csorgorevesz1981}, for $\varepsilon$ sufficiently small,
\begin{equation*}
    |B_{n}(1-\varepsilon)-B(1)|=o_{\mathbb{P}}(1) 
\end{equation*}
and
\begin{equation*}
|B_{n}(\varepsilon)-B(0)|=o_{\mathbb{P}}(1).
\end{equation*}
Which implies that
\begin{eqnarray*}
\widetilde{\mathbb{T}_{n}}(\varepsilon)&=&
\int_{U(\varepsilon)}(q^{'}(x)/q(x))B_{n}(x)dx+
o_{\mathbb{P}}(1).
\end{eqnarray*}
By \cite[Theorem 2.2.]{chengParzen1997}, we have, for $\beta>1/2$ ($\beta$ in condition (K.3)),
\begin{equation}
\sup_{u\in U(\varepsilon)}\sqrt{n}(\log\big(\widehat{q}_n(u)\big)-\log(q(u)\big)\big) = o_{\mathbb{P}}(1),
\end{equation}
which implies, by using Taylor expansion,  that
\begin{equation}
\mathbb{L}_{n}(\varepsilon)=o_{\mathbb{P}}(1),
\end{equation}
where $\mathbb{L}_{n}(\varepsilon)$ is given in (\ref{6-tilds}).
This gives
\begin{eqnarray*}
\widetilde{\mathbb{T}_{n}}(\varepsilon)+\mathbb{L}_{n}(\varepsilon)&=&\int_{U(\varepsilon)}(q(x)^{'}/q(x))B_{n}(x) dx
+o_{\mathbb{P}}(1).
\end{eqnarray*}
It follows,
\begin{eqnarray*}
\sqrt{n}(\widehat{H}_{\varepsilon;n}(X) -H_\varepsilon(X))\stackrel{d}{=}\int_{U(\varepsilon)}(q^{'}(x)/q(x))B_{n}(x)dx+o_{\mathbb{P}}(1).
\end{eqnarray*}
Here and elsewhere we denote by ``$\stackrel{d}{=}$'' the equality in distribution.
\noindent Note that $$\int_{U(\varepsilon)}(q^{'}(x)/q(x))B_{n}(x)dx$$ is a gaussian random variable with mean
\begin{eqnarray*}
\mathbb{E}\left(\int_{U(\varepsilon)}(q^{'}(u)/q(u))B_{n}(u)du\right)&=&\int_{U(\varepsilon)}(q^{'}(u)/q(u))\mathbb{E}(B_{n}(u))du\\&=&0,
\end{eqnarray*}
and variance
\begin{eqnarray*}
\lefteqn{\mathbb{E}\left(\int_{U(\varepsilon)}(q^{'}(u)/q(u))B_{n}(u)du\int_{U(\varepsilon)}(q^{'}(v)/q(v))B_{n}(v)dv\right)}\\
&=&\mathbb{E}\left(\int_{U(\varepsilon)}\int_{U(\varepsilon)}(q^{'}(u)/q(u))(q^{'}(v)/q(v))B_{n}(u)B_{n}(v)dudv\right)
\\
&=&\int_{U(\varepsilon)}\int_{U(\varepsilon)}(q^{'}(u)/q(u))(q^{'}(v)/q(v))\mathbb{E}(B_{n}(u)B_{n}(v))dudv
\\
&=&\int_{U(\varepsilon)}\int_{U(\varepsilon)}(q^{'}(u)/q(u))(q^{'}(v)/q(v))(\min(u,v)-uv)dudv
\\
&=&\int_{U(\varepsilon)}(q^{'}(v)/q(v))\left((1-v)\int_0^{v}(q^{'}(u)/q(u))udu+v\int_v^{1}(q^{'}(u)/q(u))(1-u)du\right)dv\\
&=&-\varepsilon\log(q(\varepsilon))\log(q(1-\varepsilon))+\varepsilon\log^2(q(\varepsilon))
-\log(q(1-\varepsilon))\int_{U(\varepsilon)}\log(q(x))dx\\
&&+H_\varepsilon(X)\Big(\log(q(1-\varepsilon))-H_\varepsilon(X)\Big)\\
&=&\int_{U(\varepsilon)}\log^2(q(x))dx+\varepsilon\log^2(q(\varepsilon))+\varepsilon\log^2(q(1-\varepsilon))-A^2(\varepsilon)\\
&=&\left(\int_{[0,1]}\log^2(q(x))dx+o\Big(\vartheta(\varepsilon)\Big)\right)-\left(H(X)+o\Big(\eta(\varepsilon)\Big)\right)^2\\
&=&\mbox{Var}\left\{\log\Big(q\big(F(X)\big)\Big)\right\}+o\Big(\vartheta(\varepsilon)+\eta^{2}(\varepsilon)\Big).
\end{eqnarray*}
Thus the proof is complete.\hfill$\Box$

\def\cprime{$'$}

\end{document}